\documentclass[12pt,a4paper]{amsart}

\newtheorem{theo+}              {Theorem}          % [section]
\newtheorem{prop+}  [theo+]     {Proposition}
\newtheorem{coro+}  [theo+]     {Corollary}
\newtheorem{lemm+}  [theo+]     {Lemma}
\newtheorem{exam+}  [theo+]     {Example}
\newtheorem{rema+}  [theo+]     {Remark}
\newtheorem{defi+}  [theo+]     {Definition}

\newenvironment{theorem}{\begin{theo+}}{\end{theo+}}
\newenvironment{proposition}{\begin{prop+}}{\end{prop+}}

\newenvironment{lemma}{\begin{lemm+}}{\end{lemm+}}
\newenvironment{example}{\begin{exam+}}{\end{exam+}}
\newenvironment{remark}{\begin{rema+}}{\end{rema+}}
\newenvironment{definition}{\begin{defi+}}{\end{defi+}}
\theoremstyle{definition}

\def\E{/\kern-1.0em \equiv }

\author{Ye-Lin Ou and Frederick Wilhelm$^{\ast }$}
\thanks{$^{\ast }$Support from NSF grant
DMS-0102776 is gratefully acknowledged by the second author.}
\address{Department of Mathematics\newline
\indent California University, Riverside\newline \indent
Riverside, CA 92521, U. S. A.\newline \indent
E-mail:yelino@ucr.edu (Ou),\;\;\;fred@math.ucr.edu (Wilhelm)}

\begin{document}

\title[Homothetic Submersions $\&$ Nonnegative Curvature]{Horizontally
Homothetic Submersions and Nonnegative Curvature}
\maketitle

\begin{abstract}
We show that any horizontally homothetic submersion from a compact manifold
of nonnegative sectional curvature is a Riemannian submersion.
\end{abstract}

\bigskip

The lack of examples of manifolds with positive sectional curvature has been
a major obstacle to their classification. Apart from $S^{n},$ every known
compact manifold with positive sectional curvature is constructed as the
image of a Riemannian submersion of a compact manifold with nonnegative
sectional curvature.

Here we study a generalization of Riemannian submersions called
\textquotedblleft horizontally homothetic\textquotedblright\ submersions.
For this larger class of submersions, the analog of O'Neill's horizontal
curvature equation has exactly one extra term (\cite{Gu1} and \cite{KW}).
This extra term is always nonnegative and can potentially be positive. So
the horizontal curvature equation suggests that a single horizontally
homothetic submersion is more likely to have a positively curved image than
a given Riemannian submersion. Since horizontally homothetic submersions are
(a priori) more abundant, one is lead to believe that they have much more
potential for creating positive curvature than Riemannian submersions.
Unfortunately, our main result suggests that this is an illusion.

\medskip

\noindent \textbf{Main Theorem}. \emph{Every horizontally homothetic
submersion from a compact Riemannian manifold with nonnegative sectional
curvature is a Riemannian submersion (up to a change of scale on the base
space).}\newline
\medskip

This generalizes the result in \cite{OW} that any horizontally homothetic
submersion of a round sphere with $1$--dimensional fibers is a Riemannian
submersion.

In the special case of maps from $\mathbb{R}^{n}\longrightarrow \mathbb{R},%
\mathbb{\ }$horizontally homothetic submersions appeared as solutions of the
so-called \textquotedblleft Infinity Laplace Equation\textquotedblright\
introduced by Aronsson ([Ar]) in his study of ``optimal" Lipschitz extension
of functions in the late 1960s. Nowadays these solutions are called
\textquotedblleft infinity-harmonic functions\textquotedblright , and have
been the subject of a great deal of current research (see e.g. \cite{ACJ}, 
\cite{BB}, \cite{Ba}, \cite{BEJ}, \cite{BJW1}, \cite{BJW2}, \cite{Bh}, \cite%
{CE}, \cite{CEG}, \cite{CIL}, \cite{CY}, \cite{EG}, \cite{EY}, \cite{J}, 
\cite{JK}, \cite{JLM1}, \cite{JLM2}, \cite{LM1}, \cite{LM2}, \cite{Ob}, and
the references therein) including applications in the areas of image
processing (see e.g. \cite{CMS}, \cite{Sa}), mass transfer (see e.g. \cite%
{EG}), and shape metamorphisms (see e.g. \cite{CEPB}).

In full generality, horizontally homothetic submersions arose in the study
of $p$-harmonic morphisms in [BE], [BW], [BL], [Gu1], [Gu2], [KW], [Lo],
[Ou1], [Ou2], \cite{OW}, [Ta], and [Sv]. There are many familiar examples of
horizontally homothetic submersions from incomplete manifolds with
nonnegative curvature that are not Riemannian submersions. Before
elaborating, we recall the definition from [BE] or [BW].

\begin{definition}
A submersion of Riemannian manifolds, $\pi :M\longrightarrow B,$ is called
horizontally homothetic if and only if there is a smooth function $\lambda
:M\longrightarrow \mathbb{R}$ with vertical gradient so that for all
horizontal vectors $x$ and $y$ 
\begin{equation*}
\lambda ^{2}\left\langle x,y\right\rangle _{M}=\left\langle d\pi \left(
x\right) ,d\pi \left( y\right) \right\rangle _{B}.
\end{equation*}%
$\lambda $ is called the dilation of $\pi $.

If the condition about the gradient of $\lambda $ being vertical is dropped,
then $\pi $ is called a horizontally conformal submersion. 
\end{definition}

The power of Riemannian submersions for creating positive curvature stems
from O'Neill's horizontal curvature equation, which implies that a
Riemannian submersion $\pi :M\longrightarrow N$ does not decrease the
curvature of horizontal planes. In fact, the sectional curvature of a
horizontal plane spanned by orthonormal vectors $\left\{ x,y\right\} $ in $M$
is related to the curvature of $\mathrm{span}\left\{ d\pi \left( x\right)
,d\pi \left( y\right) \right\} $ by 
\begin{equation}
\mathrm{sec}_{N}\left( d\pi \left( x\right) ,d\pi \left( y\right) \right) =%
\mathrm{sec}_{M}\left( x,y\right) +3\left\vert A_{x}y\right\vert ^{2}.
\end{equation}%
Here $A$ is O'Neill's \textquotedblleft integrability\textquotedblright\
tensor for the horizontal distribution 
\begin{equation*}
A_{x}y=\frac{1}{2}\left[ X,Y\right] ^{vert},
\end{equation*}%
where $X$ and $Y$ are arbitrary extensions of $x$ and $y$ to horizontal
vector fields. Thus if $M$ has nonnegative curvature, then $\mathrm{sec}%
_{N}\left( d\pi \left( x\right) ,d\pi \left( y\right) \right) >0$ if either 
\begin{eqnarray*}
\mathrm{sec}_{M}\left( x,y\right) &>&0\text{ or} \\
A_{x}y &\neq &0.
\end{eqnarray*}

The generalization of O'Neill's equation for horizontally conformal
submersions was discovered independently by Kasue and Washio in \cite{KW}
and Gudmundsson in \cite{Gu1}. If $\pi :M\longrightarrow N$ is horizontally
homothetic, then the sectional curvature of a horizontal plane spanned by
orthonormal vectors $\left\{ x,y\right\} $ in $M$ is related to the
curvature of $\mathrm{span}\left\{ d\pi \left( x\right) ,d\pi \left(
y\right) \right\} $ by 
\begin{equation}
\lambda ^{2}\mathrm{sec}_{N}\left( d\pi \left( x\right) ,d\pi \left(
y\right) \right) =\mathrm{sec}_{M}\left( x,y\right) +3\left\vert A\left(
x,y\right) \right\vert ^{2}+\left\vert \mathrm{grad\,ln}\lambda \right\vert
^{2}.  \label{horiz curv eqn}
\end{equation}%
(\cite{Gu1})\newline

Thus if $M$ has nonnegative curvature, then $\mathrm{sec}_{N}\left( d\pi
\left( x\right) ,d\pi \left( y\right) \right) >0$ if either 
\begin{eqnarray*}
\mathrm{sec}_{M}\left( x,y\right) &>&0, \\
\;A_{x}y &\neq &0,\text{ or} \\
\;\left\vert \mathrm{grad\,}\lambda \right\vert ^{2} &>&0.
\end{eqnarray*}%
Thus without our theorem, one would naturally suspect that horizontally
homothetic submersions have much more potential for creating positive
curvature than Riemannian submersions.

There are many familiar examples of horizontally homothetic submersions of
incomplete manifolds with nonnegative sectional curvature that are not
Riemannian submersions.

\begin{example}
Radial projection of $\mathbb{R}^{n}\setminus \left\{ 0\right\} $ onto $%
S^{n-1}$ is horizontally homothetic (see [Gu2] and also \cite{BW1} for
details).
\end{example}

\begin{example}
\label{equator}Metric projection of $S^{n}\setminus \left\{ \mathrm{north}%
\text{ \textrm{pole, south pole}}\right\} $ onto the equator is horizontally
homothetic (see [Gu2] and also \cite{BW1} for details).
\end{example}

\begin{example}
\label{join}View $S^{p+q+1}$ as the join $S^{p}\ast S^{q}.$ Then metric
projection of $S^{p+q+1}\setminus S^{p}$ onto $S^{q}$ is horizontally
homothetic.
\end{example}

For a complete example we need negative curvature.

\begin{example}
Let $H^{n}$ be hyperbolic space with the upper half space model, 
\begin{equation*}
H^{n}=\left\{ \left. \left( x_{1},\ldots ,x_{n}\right) \in \mathbb{R}%
^{n}\right\vert x_{n}>0\right\} .
\end{equation*}%
Then the projection 
\begin{equation*}
\left( x_{1},\ldots ,x_{n}\right) \longmapsto \left( x_{1},\ldots
,x_{n-1}\right)
\end{equation*}%
is a horizontally homothetic submersion onto the Euclidean space $\mathbb{R}%
^{n-1}$ (see [Gu2] and also \cite{BW1} for details).
\end{example}

\begin{example}
While the projection of a warped product onto its first factor is a
Riemannian submersion, the projection onto the second factor is horizontally
homothetic. All of the preceding examples can be viewed as the projection of
a warped product onto its second factor (see \cite{Pe}). It follows from
Theorem 1.3 in \cite{Wa} that all complete warped products of nonnegative
sectional curvature are isometric to Riemannian products. Thus our main
theorem can be viewed as a generalization of that result.
\end{example}

There also examples of horizontally homothetic submersions of incomplete
manifolds with nonnegative sectional curvature that are neither Riemannian
submersions nor warped products. For these we use projective spaces.

\begin{example}
There are multiple generalizations of Example \ref{equator} and Example \ref%
{join} that involve $\mathbb{C}P^{n},\mathbb{H}P^{n},$ and $CaP^{2}.$ For
example, metric projection of $\mathbb{C}P^{n}\setminus \left\{ pt\right\} $
onto the copy of $\mathbb{C}P^{n-1}$ at maximal distance from the point is
horizontally homothetic.
\end{example}

In the flat case we have

\begin{example}
A horizontally homothetic submersion $u:\mathbb{R}^{n}\longrightarrow 
\mathbb{R}$ satisfies the \textquotedblleft Infinity Laplace
Equation\textquotedblright\ $($see, e. g., \cite{Ar} , \cite{ACJ}, \cite{BEJ}%
, \cite{CEG}, and the references therein$)$: 
\begin{equation*}
\left\langle \mathrm{grad\,}u,\mathrm{grad\,}\left\vert \mathrm{grad\,}%
u\right\vert ^{2}\right\rangle =0.
\end{equation*}%
The solutions are called infinity-harmonic functions. Although there are
nontrivial infinity-harmonic functions on open subsets of $\mathbb{R}^{n}$,
it is not known whether all smooth globally-defined infinity-harmonic
functions on $\mathbb{R}^{n}$ are affine.
\end{example}

Horizontally homothetic submersions have played an important role in the
study of $p$-harmonic morphisms. For example, by combining results in \cite%
{BE}, \cite{BG}, \cite{BL}, and \cite{Ta1} we have\medskip\ 

\noindent \textbf{Theorem} \newline
\emph{Let }$m>n\geq 2$\emph{\ and }$\varphi :(M^{m},g)\rightarrow (N^{n},h)$%
\emph{\ be a horizontally conformal submersion.}

\begin{enumerate}
\item[(I)] \emph{If }$p=n$\emph{, then }$\varphi $\emph{\ is }$p$\emph{%
-harmonic map if and only if }$\{\varphi ^{-1}(y)\}_{y\in N}$\emph{\ is a
minimal foliation of }$(M,g)$\emph{\ of codimension }$n$\emph{.}

\item[(II)] \emph{If }$p\neq n$\emph{, then any two of the following
conditions imply the third:}

\begin{enumerate}
\item $\varphi $\emph{\ is a }$p$\emph{-harmonic map,}

\item $\{\varphi ^{-1}(y)\}_{y\in N}$\emph{\ is a minimal foliation of }$%
(M,g)$\emph{\ of codimension }$n$\emph{,}

\item $\varphi $\emph{\ is horizontally homothetic\medskip .}
\end{enumerate}
\end{enumerate}

For applications of horizontally homothetic submersions in classifying $p$%
-harmonic morphisms and biharmonic morphisms between certain model spaces
see \cite{Ou1} and \cite{Ou2}.\newline

Although horizontally homothetic submersions appear to be a much broader
class than Riemannian submersions, the presence of a nonconstant smooth
function with vertical gradient imposes a fair amount of extra structure.
For example, it is easy to show

\begin{proposition}
Let $\pi :M\longrightarrow B$ be a horizontally homothetic submersion with
dilation $\lambda$ and let $r$ be a regular value of $\lambda $ so that $%
\lambda ^{-1}\left( r\right) $ is nonempty. Then 
\begin{equation*}
\pi |_{\lambda ^{-1}\left( r\right) }:\lambda ^{-1}\left( r\right)
\longrightarrow B
\end{equation*}%
is a Riemannian submersion with respect to the intrinsic metric on $\lambda
^{-1}\left( r\right) .$
\end{proposition}

The levels of $\lambda $ also force some vanishings of the $A$--tensor of $%
\pi ,$ but before we can elaborate we must refine the definition of the $A$%
--tensor in the horizontally homothetic case.

In the Riemannian case, O'Neill defined the two fundamental tensors$,$%
\begin{eqnarray*}
A_{Z}W &=&\left( \nabla _{Z^{horiz}}W^{horiz}\right) ^{vert}+\left( \nabla
_{Z^{horiz}}W^{vert}\right) ^{horiz}\text{ and} \\
T_{Z}W &=&\left( \nabla _{Z^{vert}}W^{vert}\right) ^{horiz}+\left( \nabla
_{Z^{vert}}W^{horiz}\right) ^{vert}.
\end{eqnarray*}%
These tensors are also important for horizontally homothetic submersions,
except that it is convenient to modify the definition of the $A$--tensor to 
\begin{equation*}
A_{Z}W=\frac{1}{2}\left[ Z^{horiz},W^{horiz}\right] ^{vert}+\left( \nabla
_{Z^{horiz}}W^{vert}\right) ^{horiz,\text{ }z^{\perp }},
\end{equation*}%
where the last superscript $^{horiz,\text{ }z^{\perp }}$ indicates that we
are taking the component that is horizontal and perpendicular to $z.$ In the
Riemannian case this definition of the $A$--tensor coincides with O'Neill's
since for a Riemannian submersion $\frac{1}{2}\left[ Z^{horiz},W^{horiz}%
\right] ^{vert}=\left( \nabla _{Z^{horiz}}W^{horiz}\right) ^{vert}$ and $%
\left( \nabla _{Z^{horiz}}W^{vert}\right) ^{horiz}$ is already perpendicular
to $Z.$ On the other hand, we will show below that for horizontally
homothetic submersions, $\nabla _{Z^{horiz}}Z^{horiz}$ has a component that
is proportional to $\mathrm{grad\,}\lambda $ and (dually) $\nabla
_{Z^{horiz}}\mathrm{grad\,}\lambda $ is proportional to $Z^{horiz}.$

We assume that the reader has a working knowledge of O'Neill's foundational
paper \cite{On}. We use superscripts $^{horiz}$ and $^{vert}$ on vectors to
denote the horizontal and vertical parts. Similarly $z^{\mathrm{grad\,}%
\lambda }$ stands for the component of $z$ in the direction of $\mathrm{%
grad\,}\lambda ;$ $z^{\mathrm{grad\,}\lambda ,\perp }$ stands for the
component of $z$ that is perpendicular to $\mathrm{grad\,}\lambda ;$ and
more generally, $w^{z,\perp }$ stands for the component of $w$ that is
perpendicular to $z.$ We use 
\begin{equation*}
D_{v}\left( f\right)
\end{equation*}%
for the derivative of $f$ in the direction of $v.$

For the moment, we study an arbitrary horizontally homothetic submersion 
\begin{equation*}
\pi :M\longrightarrow B
\end{equation*}%
with dilation $\lambda .$ Later we will add the curvature and compactness
hypotheses. Let 
\begin{equation*}
\rho =-\mathrm{ln}\,\lambda .
\end{equation*}%
Then the gradient of $\rho $ is also vertical.

\begin{lemma}
Let $\gamma $ be a geodesic in $B,$ and let $X$ be the horizontal lift of $%
\dot{\gamma}$ defined on all of the horizontal lifts of $\gamma .$ Then%
\begin{equation*}
\left[ X,\frac{\mathrm{grad\,}\rho }{\left\vert \mathrm{grad\,}\rho
\right\vert ^{2}}\right] \text{ is vertical}
\end{equation*}%
and 
\begin{equation*}
\left\langle \left[ X,\frac{\mathrm{grad\,}\rho }{\left\vert \mathrm{grad\,}%
\rho \right\vert ^{2}}\right] ,\mathrm{grad\,}\rho \right\rangle =0.
\end{equation*}%
In fact, for any horizontal vector fields $X$ and $Y,$ 
\begin{equation*}
\left\langle \left[ X,Y\right] ,\mathrm{grad\,}\rho \right\rangle =0.
\end{equation*}%
Moreover, if $V$ is any vertical field that is perpendicular to $\mathrm{%
grad\,}\rho ,$ then 
\begin{equation*}
\left\langle \left[ X,V\right] ,\mathrm{grad\,}\rho \right\rangle =0.
\end{equation*}
\end{lemma}

\begin{proof}
Since $X$ is basic horizontal, 
\begin{equation*}
\left[ X,\frac{\mathrm{grad\,}\rho }{\left\vert \mathrm{grad\,}\rho
\right\vert ^{2}}\right] \text{ is vertical.}
\end{equation*}

For the second statement, note that

\begin{eqnarray*}
\left\langle \left[ X,\frac{\mathrm{grad\,}\rho }{\left\vert \mathrm{grad\,}%
\rho \right\vert ^{2}}\right] ,\mathrm{grad\,}\rho \right\rangle &=&D_{\left[
X,\frac{\mathrm{grad\,}\rho }{\left\vert \mathrm{grad\,}\rho \right\vert ^{2}%
}\right] }(\rho ) \\
&=&D_{X}(D_{\frac{\mathrm{grad\,}\rho }{\left\vert \mathrm{grad\,}\rho
\right\vert ^{2}}}(\rho ))-D_{\frac{\mathrm{grad\,}\rho }{\left\vert \mathrm{%
grad\,}\rho \right\vert ^{2}}}(D_{X}(\rho )) \\
&=&D_{X}(1)-D_{\frac{\mathrm{grad\,}\rho }{\left\vert \mathrm{grad\,}\rho
\right\vert ^{2}}}(0)=0.
\end{eqnarray*}%
In the second term of the third equality we have used the fact that $%
D_{X}(\rho )=0$ for any horizontal vector field $X,$ since the dilation $%
\lambda $ and hence $\rho $ is constant along any horizontal curve.\newline

$\left\langle \left[ X,Y\right] ,\mathrm{grad\,}\rho \right\rangle =0$ and $%
\left\langle \left[ X,V\right] ,\mathrm{grad\,}\rho \right\rangle =0$
because $X,$ $Y,$ and $V$ are tangent to submanifolds (the levels of $%
\lambda )$ that are perpendicular to $\mathrm{grad\,}\lambda .$
\end{proof}

\begin{lemma}
\label{nabla_x grad rho}For any horizontal vector $X,$ 
\begin{equation}
\nabla _{X}\frac{\mathrm{grad\,}\rho }{\left\vert \mathrm{grad\,}\rho
\right\vert ^{2}}=X-\frac{D_{X}\left\langle \mathrm{grad\,}\rho ,\mathrm{%
grad\,}\rho \right\rangle }{2\left\langle \mathrm{grad\,}\rho ,\mathrm{grad\,%
}\rho \right\rangle ^{2}}\mathrm{grad\,}\rho +\frac{1}{2}\left[ X,\frac{%
\mathrm{grad\,}\rho }{\left\vert \mathrm{grad\,}\rho \right\vert ^{2}}\right]
.  \label{AA}
\end{equation}

In particular, the $A$--tensor satisfies ${\Large \,}$ 
\begin{equation*}
A\,_{{\Huge \cdot }}\,\mathrm{grad\,}\rho =0.
\end{equation*}
\end{lemma}

\begin{proof}
To compute $\left\langle \nabla _{X}\frac{\mathrm{grad\,}\rho }{\left\vert 
\mathrm{grad\,}\rho \right\vert ^{2}},X\right\rangle $ we can now use the
Koszul formula to get 
\begin{eqnarray*}
2\left\langle \nabla _{X}\frac{\mathrm{grad\,}\rho }{\left\vert \mathrm{%
grad\,}\rho \right\vert ^{2}},X\right\rangle _{M} &=&D_{\frac{\mathrm{grad\,}%
\rho }{\left\vert \mathrm{grad\,}\rho \right\vert ^{2}}}\left\langle
X,X\right\rangle _{M} \\
&=&D_{\frac{\mathrm{grad\,}\rho }{\left\vert \mathrm{grad\,}\rho \right\vert
^{2}}}\left( e^{2\rho }\left\langle d\pi \left( X\right) ,\pi \left(
X\right) \right\rangle _{B}\right) \\
&=&\frac{\left\langle d\pi \left( X\right) ,\pi \left( X\right)
\right\rangle _{B}}{\left\vert \mathrm{grad\,}\rho \right\vert ^{2}}D_{%
\mathrm{grad\,}\rho }\left( e^{2\rho }\right) \\
&=&\frac{\left\langle d\pi \left( X\right) ,\pi \left( X\right)
\right\rangle _{B}}{\left\vert \mathrm{grad\,}\rho \right\vert ^{2}}%
2e^{2\rho }\left\vert \mathrm{grad\,}\rho \right\vert ^{2} \\
&=&2e^{2\rho }\left\langle d\pi \left( X\right) ,\pi \left( X\right)
\right\rangle _{B} \\
&=&2\left\langle X,X\right\rangle _{M}.
\end{eqnarray*}%
So 
\begin{eqnarray}
\left\langle \nabla _{X}\frac{\mathrm{grad\,}\rho }{\left\vert \mathrm{grad\,%
}\rho \right\vert ^{2}},\frac{X}{\left\vert X\right\vert }\right\rangle 
\frac{X}{\left\vert X\right\vert } &=&\left\langle \nabla _{X}\frac{\mathrm{%
grad\,}\rho }{\left\vert \mathrm{grad\,}\rho \right\vert ^{2}}%
,X\right\rangle \frac{X}{\left\vert X\right\vert ^{2}}  \label{AB} \\
&=&\left\langle X,X\right\rangle _{M}\frac{X}{\left\vert X\right\vert ^{2}}%
=X.  \notag
\end{eqnarray}%
Similarly, if $Z$ is any basic horizontal field that is perpendicular to $X$%
, then the previous Lemma and the Koszul formula give 
\begin{equation*}
2\left\langle \nabla _{X}\frac{\mathrm{grad\,}\rho }{\left\vert \mathrm{%
grad\,}\rho \right\vert ^{2}},Z\right\rangle =0.
\end{equation*}%
If $V$ is any vertical field that is perpendicular to $\frac{\mathrm{grad\,}%
\rho }{\left\vert \mathrm{grad\,}\rho \right\vert ^{2}}$ and $\left[ X,\frac{%
\mathrm{grad\,}\rho }{\left\vert \mathrm{grad\,}\rho \right\vert ^{2}}\right]
$, then the Koszul formula gives%
\begin{equation*}
2\left\langle \nabla _{X}\frac{\mathrm{grad\,}\rho }{\left\vert \mathrm{%
grad\,}\rho \right\vert ^{2}},V\right\rangle =0.
\end{equation*}

In the direction of $\frac{\mathrm{grad\,}\rho }{\left\vert \mathrm{grad\,}%
\rho \right\vert ^{2}},$ we have%
\begin{eqnarray*}
2\left\langle \nabla _{X}\frac{\mathrm{grad\,}\rho }{\left\vert \mathrm{%
grad\,}\rho \right\vert ^{2}},\frac{\mathrm{grad\,}\rho }{\left\vert \mathrm{%
grad\,}\rho \right\vert ^{2}}\right\rangle &=&D_{X}\frac{\left\langle 
\mathrm{grad\,}\rho ,\mathrm{grad\,}\rho \right\rangle }{\left\langle 
\mathrm{grad\,}\rho ,\mathrm{grad\,}\rho \right\rangle ^{2}}%
=D_{X}\left\langle \mathrm{grad\,}\rho ,\mathrm{grad\,}\rho \right\rangle
^{-1} \\
&=&-\left\langle \mathrm{grad\,}\rho ,\mathrm{grad\,}\rho \right\rangle
^{-2}D_{X}\left\langle \mathrm{grad\,}\rho ,\mathrm{grad\,}\rho \right\rangle
\\
&=&-\frac{D_{X}\left\langle \mathrm{grad\,}\rho ,\mathrm{grad\,}\rho
\right\rangle }{\left\langle \mathrm{grad\,}\rho ,\mathrm{grad\,}\rho
\right\rangle ^{2}}.
\end{eqnarray*}%
So 
\begin{equation}
\left\langle \nabla _{X}\frac{\mathrm{grad\,}\rho }{\left\vert \mathrm{grad\,%
}\rho \right\vert ^{2}},\frac{\mathrm{grad\,}\rho }{\left\vert \mathrm{grad\,%
}\rho \right\vert }\right\rangle \frac{\mathrm{grad\,}\rho }{\left\vert 
\mathrm{grad\,}\rho \right\vert }=-\frac{D_{X}\left\langle \mathrm{grad\,}%
\rho ,\mathrm{grad\,}\rho \right\rangle }{2\left\langle \mathrm{grad\,}\rho ,%
\mathrm{grad\,}\rho \right\rangle ^{2}}\mathrm{grad\,}\rho .  \label{AC}
\end{equation}%
In case $\left[ X,\frac{\mathrm{grad\,}\rho }{\left\vert \mathrm{grad\,}\rho
\right\vert ^{2}}\right] =0,$ we are done. If not, take $W=\frac{\left[ X,%
\frac{\mathrm{grad\,}\rho }{\left\vert \mathrm{grad\,}\rho \right\vert ^{2}}%
\right] }{\left\vert \left[ X,\frac{\mathrm{grad\,}\rho }{\left\vert \mathrm{%
grad\,}\rho \right\vert ^{2}}\right] \right\vert }.$ Then 
\begin{equation}
2\left\langle \nabla _{X}\frac{\mathrm{grad\,}\rho }{\left\vert \mathrm{%
grad\,}\rho \right\vert ^{2}},W\right\rangle =\left\langle \left[ X,\frac{%
\mathrm{grad\,}\rho }{\left\vert \mathrm{grad\,}\rho \right\vert ^{2}}\right]
,W\right\rangle .  \label{AD}
\end{equation}

From (\ref{AB}), (\ref{AC}) and (\ref{AD}) we obtain Equation (\ref{AA}).%
\newline

The statement about the $A$--tensor follows from its (new) definition and
the fact that $\left( \nabla _{X}\frac{\mathrm{grad\,}\rho }{\left\vert 
\mathrm{grad\,}\rho \right\vert ^{2}}\right) ^{horiz}=X.$
\end{proof}

Since $\rho$ is constant along any horizontal geodesic we have

\begin{proposition}
If the global maximum or minimum of $\rho $ occurs at $p$, then it also
occurs along the entire length of any horizontal curve passing through $p.$
\end{proposition}

In fact, we have

\begin{lemma}
\label{reg horiz}If $r$ is a regular point of $\rho $, then so is the entire
length of any horizontal curve passing through $r.$
\end{lemma}

\begin{proof}
Let $\gamma $ be a curve in the base. Then the horizontal lifts of $\gamma $
give a family of diffeomorphisms 
\begin{equation*}
H_{t}:\pi ^{-1}\left( \gamma \left( 0\right) \right) \longrightarrow \pi
^{-1}\left( \gamma \left( t\right) \right) ,\;\;H_{t}(x)={\tilde{\gamma}}(1),
\end{equation*}%
where ${\tilde{\gamma}}$ is the horizontal lift of $\gamma $ starting at $%
x\in \pi ^{-1}(\gamma (0))$. These were studied as early as 1960 in [He] and
were called \textquotedblleft Holonomy Displacement\textquotedblright\ maps
in [GG].\newline
Since a horizontal curve stays in a fixed level of $\rho ,$ it follows that
for any $x\in \pi ^{-1}\left( \gamma \left( 0\right) \right) ,$ $\rho \circ
H_{t}\left( x\right) $ is independent of $t$. Thus for a regular point $r$
of $\rho $ 
\begin{equation*}
d\rho \left( \left( dH_{t}\right) _{r}\left( \mathrm{grad\,}\rho \right)
\right)
\end{equation*}%
is independent of $t.$ Since $\mathrm{grad\,}\rho $ is nonzero at $r$, it
follows that 
\begin{equation*}
d\rho \left( \left( dH_{t}\right) _{r}\left( \mathrm{grad\,}\rho \right)
\right) \neq 0
\end{equation*}%
for all $t.$ In particular, $\mathrm{grad\,}\rho $ is nonzero for all $t.$
\end{proof}

Generalizing Lemma 2.2 in \cite{OW} we have

\begin{proposition}
Let $\tilde{\gamma}$ be a horizontal lift of a geodesic $\gamma $ in $B.$
Then, $\tilde{\gamma}$ is an intrinsic geodesic of its level set of $\rho ,$
and a geodesic of $M$ only if the level is critical.
\end{proposition}

\begin{proof}
Let $X$ be the horizontal lift of $\gamma ^{\prime }$, and let $Z$ be any
vertical field that is tangent to a level set of $\rho .$ Then using the
Koszul formula 
\begin{equation*}
2\left\langle \nabla _{X}X,Z\right\rangle =-D_{Z}\left\langle
X,X\right\rangle -\left\langle \left[ X,Z\right] ,X\right\rangle
+\left\langle \left[ Z,X\right] ,X\right\rangle .
\end{equation*}%
All three of these terms are $0.$ Indeed $D_{Z}\left\langle X,X\right\rangle
=0,$ since $Z$ is tangent to a level of $\rho $ and $X$ is basic horizontal
and $\left\langle \left[ X,Z\right] ,X\right\rangle =0,$ since $X$ is basic
horizontal.

If $Z$ is horizontal, then by the Koszul formula 
\begin{equation*}
2\left\langle \nabla _{X}X,Z\right\rangle =2D_{X}\left\langle
X,Z\right\rangle -D_{Z}\left\langle X,X\right\rangle -\left\langle \left[ X,Z%
\right] ,X\right\rangle +\left\langle \left[ Z,X\right] ,X\right\rangle ,
\end{equation*}%
and the righthand side is $0,$ since it is equal to a multiple of 
\begin{equation*}
2\left\langle \nabla _{d\pi \left( X\right) }d\pi \left( X\right) ,d\pi
\left( Z\right) \right\rangle ,
\end{equation*}%
which is $0,$ since $X$ is a lift of a geodesic field.

It follows that $\tilde{\gamma}$ is an intrinsic geodesic for its level set
of $\rho .$

On the other hand, since $X$ is basic horizontal 
\begin{eqnarray*}
\left\langle \nabla _{X}X,\frac{\mathrm{grad\,}\rho }{\left\vert \mathrm{%
grad\,}\rho \right\vert ^{2}}\right\rangle &=&-\left\langle X,\nabla _{X}%
\frac{\mathrm{grad\,}\rho }{\left\vert \mathrm{grad\,}\rho \right\vert ^{2}}%
\right\rangle \\
&=&-\left\langle X,\nabla _{\frac{\mathrm{grad\,}\rho }{\left\vert \mathrm{%
grad\,}\rho \right\vert ^{2}}}X\right\rangle \\
&=&-\frac{1}{2}D_{\frac{\mathrm{grad\,}\rho }{\left\vert \mathrm{grad\,}\rho
\right\vert ^{2}}}\left\langle X,X\right\rangle \\
&=&-\frac{1}{2}D_{\frac{\mathrm{grad\,}\rho }{\left\vert \mathrm{grad\,}\rho
\right\vert ^{2}}}\left( e^{2\rho }\left\langle d\pi \left( X\right) ,d\pi
\left( X\right) \right\rangle _{B}\right) \\
&=&-\frac{2\left\langle d\pi \left( X\right) ,d\pi \left( X\right)
\right\rangle _{B}}{2\left\vert \mathrm{grad\,}\rho \right\vert ^{2}}%
e^{2\rho }D_{\mathrm{grad\,}\rho }\left( \rho \right) \\
&=&-e^{2\rho }\left\langle d\pi \left( X\right) ,d\pi \left( X\right)
\right\rangle _{B} \\
&=&-\left\langle X,X\right\rangle _{M}
\end{eqnarray*}%
and hence $\nabla _{X}X$ is nonzero wherever $\mathrm{grad\,}\rho $ is
nonzero.
\end{proof}

\begin{lemma}
\begin{equation*}
T_{\frac{\mathrm{grad\,}\rho }{\left\vert \mathrm{grad\,}\rho \right\vert }%
}X=-\frac{D_{X}\left\langle \mathrm{grad\,}\rho ,\mathrm{grad\,}\rho
\right\rangle }{2\left\langle \mathrm{grad\,}\rho ,\mathrm{grad\,}\rho
\right\rangle }\left( \frac{\mathrm{grad\,}\rho }{\left\vert \mathrm{grad\,}%
\rho \right\vert }\right) +\frac{1}{2}\left\vert \mathrm{grad\,}\rho
\right\vert \left[ \frac{\mathrm{grad\,}\rho }{\left\vert \mathrm{grad\,}%
\rho \right\vert ^{2}},X\right] .
\end{equation*}
\end{lemma}

\begin{proof}
Since $\left[ \frac{\mathrm{grad\,}\rho }{\left\vert \mathrm{grad\,}\rho
\right\vert ^{2}},X\right] $ is perpendicular to $\mathrm{grad\,}\rho ,$ we
use Lemma \ref{nabla_x grad rho} to get%
\begin{eqnarray*}
\left\langle T_{\frac{\mathrm{grad\,}\rho }{\left\vert \mathrm{grad\,}\rho
\right\vert }}X,\frac{\mathrm{grad\,}\rho }{\left\vert \mathrm{grad\,}\rho
\right\vert }\right\rangle &=&\left\langle \nabla _{X}\frac{\mathrm{grad\,}%
\rho }{\left\vert \mathrm{grad\,}\rho \right\vert ^{2}},\mathrm{grad\,}\rho
\right\rangle \\
&=&\left\langle X-\frac{D_{X}\left\langle \mathrm{grad\,}\rho ,\mathrm{grad\,%
}\rho \right\rangle }{2\left\langle \mathrm{grad\,}\rho ,\mathrm{grad\,}\rho
\right\rangle ^{2}}\mathrm{grad\,}\rho +\frac{1}{2}\left[ X,\frac{\mathrm{%
grad\,}\rho }{\left\vert \mathrm{grad\,}\rho \right\vert ^{2}}\right] ,%
\mathrm{grad\,}\rho \right\rangle \\
&=&-\frac{D_{X}\left\langle \mathrm{grad\,}\rho ,\mathrm{grad\,}\rho
\right\rangle }{2\left\langle \mathrm{grad\,}\rho ,\mathrm{grad\,}\rho
\right\rangle ^{2}}\left\langle \mathrm{grad\,}\rho ,\mathrm{grad\,}\rho
\right\rangle \\
&=&-\frac{D_{X}\left\langle \mathrm{grad\,}\rho ,\mathrm{grad\,}\rho
\right\rangle }{2\left\langle \mathrm{grad\,}\rho ,\mathrm{grad\,}\rho
\right\rangle }.
\end{eqnarray*}%
If $V$ is a vertical vector that is perpendicular to $\mathrm{grad\,}\rho $
and $\left[ \frac{\mathrm{grad\,}\rho }{\left\vert \mathrm{grad\,}\rho
\right\vert ^{2}},X\right] ,$ then the Koszul formula gives%
\begin{equation*}
\left\langle T_{\frac{\mathrm{grad\,}\rho }{\left\vert \mathrm{grad\,}\rho
\right\vert }}X,V\right\rangle =\left\langle \nabla _{\frac{\mathrm{grad\,}%
\rho }{\left\vert \mathrm{grad\,}\rho \right\vert }}X,V\right\rangle =0.
\end{equation*}%
Finally, if $W=\frac{\left[ \frac{\mathrm{grad\,}\rho }{\left\vert \mathrm{%
grad\,}\rho \right\vert ^{2}},X\right] }{\left\vert \left[ \frac{\mathrm{%
grad\,}\rho }{\left\vert \mathrm{grad\,}\rho \right\vert ^{2}},X\right]
\right\vert },$ then%
\begin{eqnarray*}
2\left\langle T_{\frac{\mathrm{grad\,}\rho }{\left\vert \mathrm{grad\,}\rho
\right\vert }}X,W\right\rangle &=&2\left\vert \mathrm{grad\,}\rho
\right\vert \left\langle \nabla _{\frac{\mathrm{grad\,}\rho }{\left\vert 
\mathrm{grad\,}\rho \right\vert ^{2}}}X,W\right\rangle \\
&=&\left\vert \mathrm{grad\,}\rho \right\vert \left\langle \left[ \frac{%
\mathrm{grad\,}\rho }{\left\vert \mathrm{grad\,}\rho \right\vert ^{2}},X%
\right] ,W\right\rangle .
\end{eqnarray*}
\end{proof}

The following Lemma can be found in \cite{BW1}.

\begin{lemma}
For any basic horizontal field $X$ and any vertical field $V,$%
\begin{eqnarray*}
\left\langle R\left( X,V\right) V,X\right\rangle &=&\left\langle \left(
\nabla _{X}T\right) _{V}V,X\right\rangle -\left\langle
T_{V}X,T_{V}X\right\rangle +\left\langle A_{X}V,A_{X}V\right\rangle \\
&&-\left\langle \nabla _{V}\left( \nabla _{X}V^{\mathrm{grad\,}\rho }\right)
^{horiz},X\right\rangle +\left\langle \nabla _{X}\left( \nabla _{V}V\right)
^{vert},X\right\rangle .
\end{eqnarray*}
\end{lemma}

\begin{remark}
Note that the first three terms are precisely O'Neill's formula for
vertizontal curvature of a Riemannian submersion. In that case, $%
\left\langle \nabla _{V}\left( \nabla _{X}V^{\mathrm{grad\,}\rho }\right)
^{horiz},X\right\rangle $ vanishes since $\mathrm{grad\,}\rho =0,$ and the
last term vanishes by the antisymmetry of the $A$--tensor.
\end{remark}

\begin{proof}
\begin{eqnarray*}
\left\langle R\left( X,V\right) V,X\right\rangle &=&\left\langle \nabla
_{X}\nabla _{V}V,X\right\rangle -\left\langle \nabla _{V}\nabla
_{X}V,X\right\rangle -\left\langle \nabla _{\left[ X,V\right]
}V,X\right\rangle \\
&=&\left\langle \nabla _{X}\left( \nabla _{V}V\right)
^{horiz},X\right\rangle +\left\langle \nabla _{X}\left( \nabla _{V}V\right)
^{vert},X\right\rangle \\
&&-\left\langle \nabla _{V}\left( \nabla _{X}V\right)
^{horiz},X\right\rangle -\left\langle \nabla _{V}\left( \nabla _{X}V\right)
^{vert},X\right\rangle \\
&&-\left\langle \nabla _{\left[ X,V\right] ^{vert}}V,X\right\rangle .
\end{eqnarray*}%
Letting the superscript $^{\mathrm{grad\,}\rho }$ denote the component
tangent to $\mathrm{grad\,}\rho ,$ and $^{\mathrm{grad\,}\rho ,\perp }$
denote the component perpendicular to $\mathrm{grad\,}\rho ,$ and using the
definition of the $A$--tensor and Lemma \ref{nabla_x grad rho}, we get that
the third term is 
\begin{eqnarray*}
-\left\langle \nabla _{V}\left( \nabla _{X}V\right) ^{horiz},X\right\rangle
&=&-\left\langle \nabla _{V}\left( \nabla _{X}V^{\mathrm{grad\,}\rho
}\right) ^{horiz},X\right\rangle -\left\langle \nabla _{V}\left( \nabla
_{X}V^{\mathrm{grad\,}\rho ,\perp }\right) ^{horiz},X\right\rangle \\
&=&-\left\langle \nabla _{V}\left( \nabla _{X}V^{\mathrm{grad\,}\rho
}\right) ^{horiz},X\right\rangle -\left\langle \nabla _{V}\left(
A_{X}V\right) ,X\right\rangle .
\end{eqnarray*}

Thus if we also use the definition of the $T$--tensor, we get 
\begin{eqnarray*}
\left\langle R\left( X,V\right) V,X\right\rangle &=&\left\langle \nabla
_{X}\left( T_{V}V\right) ,X\right\rangle -\left\langle \nabla _{V}\left(
A_{X}V\right) ,X\right\rangle \\
&&-\left\langle T_{V}\left( \nabla _{X}V\right) ^{vert},X\right\rangle
-\left\langle \nabla _{\left( \nabla _{X}V\right) ^{vert}}V,X\right\rangle
+\left\langle \nabla _{\left( \nabla _{V}X\right) ^{vert}}V,X\right\rangle \\
&&-\left\langle \nabla _{V}\left( \nabla _{X}V^{\mathrm{grad\,}\rho }\right)
^{horiz},X\right\rangle +\left\langle \nabla _{X}\left( \nabla _{V}V\right)
^{vert},X\right\rangle \\
&=&\left\langle \nabla _{X}\left( T_{V}V\right) ,X\right\rangle
-\left\langle T_{V}\left( \nabla _{X}V\right) ^{vert},X\right\rangle
-\left\langle T_{\left( \nabla _{X}V\right) ^{vert}}V,X\right\rangle \\
&&+\left\langle T_{\left( T_{V}X\right) ^{vert}}V,X\right\rangle
+\left\langle A_{X}V,A_{X}V\right\rangle \\
&&-\left\langle \nabla _{V}\left( \nabla _{X}V^{\mathrm{grad\,}\rho }\right)
^{horiz},X\right\rangle +\left\langle \nabla _{X}\left( \nabla _{V}V\right)
^{vert},X\right\rangle \\
&=&\left\langle \left( \nabla _{X}T\right) _{V}V,X\right\rangle
-\left\langle T_{V}X,T_{V}X\right\rangle +\left\langle
A_{X}V,A_{X}V\right\rangle \\
&&-\left\langle \nabla _{V}\left( \nabla _{X}V^{\mathrm{grad\,}\rho }\right)
^{horiz},X\right\rangle +\left\langle \nabla _{X}\left( \nabla _{V}V\right)
^{vert},X\right\rangle .
\end{eqnarray*}
\end{proof}

When $V=\frac{\mathrm{grad\,}\rho }{\left\vert \mathrm{grad\,}\rho
\right\vert }$, the first non-Riemannian term can be simplified as follows.

\begin{proposition}
When $V=\frac{\mathrm{grad\,}\rho }{\left\vert \mathrm{grad\,}\rho
\right\vert },$ 
\begin{equation*}
\left\langle \nabla _{V}\left( \nabla _{X}V^{\mathrm{grad\,}\rho }\right)
^{horiz},X\right\rangle =\left( D_{\frac{\mathrm{grad\,}\rho }{\left\vert 
\mathrm{grad\,}\rho \right\vert }}\left\vert \mathrm{grad\,}\rho \right\vert
+\left\vert \mathrm{grad\,}\rho \right\vert ^{2}\right) \left\langle
X,X\right\rangle
\end{equation*}
\end{proposition}

\begin{proof}
We use Equation (\ref{AA}) to get 
\begin{eqnarray*}
\left\langle \nabla _{\frac{\mathrm{grad\,}\rho }{\left\vert \mathrm{grad\,}%
\rho \right\vert }}\left( \nabla _{X}\frac{\mathrm{grad\,}\rho }{\left\vert 
\mathrm{grad\,}\rho \right\vert }\right) ^{horiz},X\right\rangle
&=&\left\langle \nabla _{\frac{\mathrm{grad\,}\rho }{\left\vert \mathrm{%
grad\,}\rho \right\vert }}\left( \left\vert \mathrm{grad\,}\rho \right\vert
\nabla _{X}\frac{\mathrm{grad\,}\rho }{\left\vert \mathrm{grad\,}\rho
\right\vert ^{2}}\right) ^{horiz},X\right\rangle \\
&=&\left\langle \nabla _{\frac{\mathrm{grad\,}\rho }{\left\vert \mathrm{%
grad\,}\rho \right\vert }}\left( \left\vert \mathrm{grad\,}\rho \right\vert
X\right) ^{horiz},X\right\rangle \\
&=&\left\langle \left(D_{\frac{\mathrm{grad\,}\rho }{\left\vert \mathrm{%
grad\,}\rho \right\vert }} \left\vert \mathrm{grad\,}\rho \right\vert
\right) X+\left\vert \mathrm{grad\,}\rho \right\vert \nabla _{\frac{\mathrm{%
grad\,}\rho }{\left\vert \mathrm{grad\,}\rho \right\vert }}X,X\right\rangle
\\
&=&\left(D_{\frac{\mathrm{grad\,}\rho }{\left\vert \mathrm{grad\,}\rho
\right\vert }} \left\vert \mathrm{grad\,}\rho \right\vert \right)
\left\langle X,X\right\rangle +\left\vert \mathrm{grad\,}\rho \right\vert
^{2}\left\langle \nabla _{\frac{\mathrm{grad\,}\rho }{\left\vert \mathrm{%
grad\,}\rho \right\vert ^{2}}}X,X\right\rangle \\
&=&\left(D_{\frac{\mathrm{grad\,}\rho }{\left\vert \mathrm{grad\,}\rho
\right\vert }} \left\vert \mathrm{grad\,}\rho \right\vert \right)
\left\langle X,X\right\rangle +\left\vert \mathrm{grad\,}\rho \right\vert
^{2}\left\langle X,X\right\rangle .
\end{eqnarray*}
\end{proof}

The last non-Riemannian curvature term can also be simplified when $V=\frac{%
\mathrm{grad\,}\rho }{\left\vert \mathrm{grad\,}\rho \right\vert }.$

\begin{proposition}
When $V=$\bigskip $\frac{\mathrm{grad\,}\rho }{\left\vert \mathrm{grad\,}%
\rho \right\vert },$%
\begin{equation*}
\left\langle \nabla _{X}\left( \nabla _{V}V\right) ^{vert},X\right\rangle =0.
\end{equation*}
\end{proposition}

\begin{proof}
We have

\begin{eqnarray}
\left\langle \nabla _{X}\left( \nabla _{V}V\right) ^{vert},X\right\rangle
&=&\left\langle \nabla _{X}\left( \nabla _{\frac{\mathrm{grad\,}\rho }{%
\left\vert \mathrm{grad\,}\rho \right\vert }}\frac{\mathrm{grad\,}\rho }{%
\left\vert \mathrm{grad\,}\rho \right\vert }\right) ^{vert,\perp
},X\right\rangle +  \notag \\
&&\left\langle \nabla _{X}\left( \nabla _{\frac{\mathrm{grad\,}\rho }{%
\left\vert \mathrm{grad\,}\rho \right\vert }}\frac{\mathrm{grad\,}\rho }{%
\left\vert \mathrm{grad\,}\rho \right\vert }\right) ^{\mathrm{grad\,}\rho
},X\right\rangle .  \notag
\end{eqnarray}%
The first term on the right-hand side of the above equation is $0$ by the
antisymmetry of the $A$--tensor. Thus, 
\begin{equation*}
\left\langle \nabla _{X}\left( \nabla _{V}V\right) ^{vert},X\right\rangle
=\left\langle \nabla _{X}\left( \nabla _{\frac{\mathrm{grad\,}\rho }{%
\left\vert \mathrm{grad\,}\rho \right\vert }}\frac{\mathrm{grad\,}\rho }{%
\left\vert \mathrm{grad\,}\rho \right\vert }\right) ^{\mathrm{grad\,}\rho
},X\right\rangle .
\end{equation*}%
Since 
\begin{equation*}
\left\langle \nabla _{\frac{\mathrm{grad\,}\rho }{\left\vert \mathrm{grad\,}%
\rho \right\vert }}\frac{\mathrm{grad\,}\rho }{\left\vert \mathrm{grad\,}%
\rho \right\vert },\frac{\mathrm{grad\,}\rho }{\left\vert \mathrm{grad\,}%
\rho \right\vert }\right\rangle =\frac{1}{2}D_{\frac{\mathrm{grad\,}\rho }{%
\left\vert \mathrm{grad\,}\rho \right\vert }}\left\langle \frac{\mathrm{%
grad\,}\rho }{\left\vert \mathrm{grad\,}\rho \right\vert },\frac{\mathrm{%
grad\,}\rho }{\left\vert \mathrm{grad\,}\rho \right\vert }\right\rangle =0,
\end{equation*}%
the entire last term is $0.$
\end{proof}

Thus, when $V=\frac{\mathrm{grad\,}\rho }{\left\vert \mathrm{grad\,}\rho
\right\vert }$ \bigskip we get 
\begin{eqnarray}
&&\left\langle R\left( X,\frac{\mathrm{grad\,}\rho }{\left\vert \mathrm{%
grad\,}\rho \right\vert }\right) \frac{\mathrm{grad\,}\rho }{\left\vert 
\mathrm{grad\,}\rho \right\vert },X\right\rangle  \notag \\
&=&\left\langle \left( \nabla _{X}T\right) _{\frac{\mathrm{grad\,}\rho }{%
\left\vert \mathrm{grad\,}\rho \right\vert }}\frac{\mathrm{grad\,}\rho }{%
\left\vert \mathrm{grad\,}\rho \right\vert },X\right\rangle -\left\langle T_{%
\frac{\mathrm{grad\,}\rho }{\left\vert \mathrm{grad\,}\rho \right\vert }%
}X,T_{\frac{\mathrm{grad\,}\rho }{\left\vert \mathrm{grad\,}\rho \right\vert 
}}X\right\rangle +  \notag \\
&&\left\langle A_{X}\frac{\mathrm{grad\,}\rho }{\left\vert \mathrm{grad\,}%
\rho \right\vert },A_{X}\frac{\mathrm{grad\,}\rho }{\left\vert \mathrm{grad\,%
}\rho \right\vert }\right\rangle -\left( D_{\frac{\mathrm{grad\,}\rho }{%
\left\vert \mathrm{grad\,}\rho \right\vert }}\left\vert \mathrm{grad\,}\rho
\right\vert +\left\vert \mathrm{grad\,}\rho \right\vert ^{2}\right)
\left\langle X,X\right\rangle .  \notag
\end{eqnarray}%
Since the $A$--tensor term vanishes, we get

\begin{eqnarray}
&&\left\langle R\left( X,\frac{\mathrm{grad\,}\rho }{\left\vert \mathrm{%
grad\,}\rho \right\vert }\right) \frac{\mathrm{grad\,}\rho }{\left\vert 
\mathrm{grad\,}\rho \right\vert },X\right\rangle  \label{AF} \\
&=&\left\langle \nabla _{X}\left( T_{\frac{\mathrm{grad\,}\rho }{\left\vert 
\mathrm{grad\,}\rho \right\vert }}\frac{\mathrm{grad\,}\rho }{\left\vert 
\mathrm{grad\,}\rho \right\vert }\right) ,X\right\rangle -2\left\langle T_{%
\frac{\mathrm{grad\,}\rho }{\left\vert \mathrm{grad\,}\rho \right\vert }%
}\left( \nabla _{X}\frac{\mathrm{grad\,}\rho }{\left\vert \mathrm{grad\,}%
\rho \right\vert }\right) ^{vert},X\right\rangle -  \notag \\
&&-\left\langle T_{\frac{\mathrm{grad\,}\rho }{\left\vert \mathrm{grad\,}%
\rho \right\vert }}X,T_{\frac{\mathrm{grad\,}\rho }{\left\vert \mathrm{grad\,%
}\rho \right\vert }}X\right\rangle -\left( D_{\frac{\mathrm{grad\,}\rho }{%
\left\vert \mathrm{grad\,}\rho \right\vert }}\left\vert \mathrm{grad\,}\rho
\right\vert +\left\vert \mathrm{grad\,}\rho \right\vert ^{2}\right)
\left\langle X,X\right\rangle .  \notag
\end{eqnarray}

The second term on the right hand side of the previous display is given by
the following result.

\begin{proposition}
\label{P20} 
\begin{equation*}
-2\left\langle T_{\frac{\mathrm{grad\,}\rho }{\left\vert \mathrm{grad\,}\rho
\right\vert }}\left( \nabla _{X}\frac{\mathrm{grad\,}\rho }{\left\vert 
\mathrm{grad\,}\rho \right\vert }\right) ^{vert},X\right\rangle =-\frac{1}{2}%
\left\vert \mathrm{grad\,}\rho \right\vert ^{2}\left\vert \left[ X,\frac{%
\mathrm{grad\,}\rho }{\left\vert \mathrm{grad\,}\rho \right\vert ^{2}}\right]
\right\vert ^{2}.
\end{equation*}
\end{proposition}

\begin{proof}
First note that since $\frac{\mathrm{grad\,}\rho }{\left\vert \mathrm{grad\,}%
\rho \right\vert }$ has constant length, $\left\langle \nabla _{X}\frac{%
\mathrm{grad\,}\rho }{\left\vert \mathrm{grad\,}\rho \right\vert },\frac{%
\mathrm{grad\,}\rho }{\left\vert \mathrm{grad\,}\rho \right\vert }%
\right\rangle =0$. Thus 
\begin{equation*}
\left( \nabla _{X}\frac{\mathrm{grad\,}\rho }{\left\vert \mathrm{grad\,}\rho
\right\vert }\right) ^{vert}=\left\vert \mathrm{grad\,}\rho \right\vert
\left( \nabla _{X}\frac{\mathrm{grad\,}\rho }{\left\vert \mathrm{grad\,}\rho
\right\vert ^{2}}\right) ^{vert,\perp },
\end{equation*}%
where the superscript $^{vert,\perp }$ signifies the component that is
vertical and perpendicular to $\mathrm{grad\,}\rho $. Using Equation (\ref%
{AA}) we obtain 
\begin{equation*}
\left( \nabla _{X}\frac{\mathrm{grad\,}\rho }{\left\vert \mathrm{grad\,}\rho
\right\vert }\right) ^{vert}=\left\vert \mathrm{grad\,}\rho \right\vert 
\frac{1}{2}\left[ X,\frac{\mathrm{grad\,}\rho }{\left\vert \mathrm{grad\,}%
\rho \right\vert ^{2}}\right] .
\end{equation*}%
Therefore 
\begin{equation*}
\left\langle T_{\frac{\mathrm{grad\,}\rho }{\left\vert \mathrm{grad\,}\rho
\right\vert }}\left( \nabla _{X}\frac{\mathrm{grad\,}\rho }{\left\vert 
\mathrm{grad\,}\rho \right\vert }\right) ^{vert},X\right\rangle =\frac{1}{2}%
\left\langle T_{\mathrm{grad\,}\rho }\left[ X,\frac{\mathrm{grad\,}\rho }{%
\left\vert \mathrm{grad\,}\rho \right\vert ^{2}}\right] ,X\right\rangle .
\end{equation*}%
If we set $W=\left[ X,\frac{\mathrm{grad\,}\rho }{\left\vert \mathrm{grad\,}%
\rho \right\vert ^{2}}\right] $, then we get 
\begin{eqnarray*}
\left\langle T_{\frac{\mathrm{grad\,}\rho }{\left\vert \mathrm{grad\,}\rho
\right\vert }}\left( \nabla _{X}\frac{\mathrm{grad\,}\rho }{\left\vert 
\mathrm{grad\,}\rho \right\vert }\right) ^{vert},X\right\rangle &=&\frac{1}{2%
}\left\langle T_{\mathrm{grad\,}\rho }W,X\right\rangle \\
&=&\frac{1}{2}\left\langle \nabla _{\mathrm{grad\,}\rho }W,X\right\rangle \\
&=&\frac{1}{4}\left\langle \left[ X,\mathrm{grad\,}\rho \right]
,W\right\rangle \\
&=&\frac{1}{4}\left\langle \left[ X,\mathrm{grad\,}\rho \right] ,\left[ X,%
\frac{\mathrm{grad\,}\rho }{\left\vert \mathrm{grad\,}\rho \right\vert ^{2}}%
\right] \right\rangle .
\end{eqnarray*}%
Since $\left[ X,\frac{\mathrm{grad\,}\rho }{\left\vert \mathrm{grad\,}\rho
\right\vert ^{2}}\right] $ is perpendicular to $\mathrm{grad\,}\rho $ we get 
\begin{equation*}
\left\langle T_{\frac{\mathrm{grad\,}\rho }{\left\vert \mathrm{grad\,}\rho
\right\vert }}\left( \nabla _{X}\frac{\mathrm{grad\,}\rho }{\left\vert 
\mathrm{grad\,}\rho \right\vert }\right) ^{vert},X\right\rangle =\frac{1}{4}%
\left\vert \mathrm{grad\,}\rho \right\vert ^{2}\left\langle \left[ X,\frac{%
\mathrm{grad\,}\rho }{\left\vert \mathrm{grad\,}\rho \right\vert ^{2}}\right]
,\left[ X,\frac{\mathrm{grad\,}\rho }{\left\vert \mathrm{grad\,}\rho
\right\vert ^{2}}\right] \right\rangle ,
\end{equation*}%
as desired.
\end{proof}

Using Lemma \ref{reg horiz}, continuity of the second derivatives of $\rho ,$
and compactness give us

\begin{proposition}
\label{P21} Let $\pi :M\longrightarrow B$ be a horizontally homothetic
submersion of a compact Riemannian manifold with nonconstant dilation
function.

Throughout any horizontal lift of a geodesic in $B$ that passes through
regular points of $\rho $ and is sufficiently close to the minimum level of $%
\rho ,$%
\begin{equation*}
-\left( D_{\frac{\mathrm{grad\,}\rho }{\left\vert \mathrm{grad\,}\rho
\right\vert }}\left\vert \mathrm{grad\,}\rho \right\vert +\left\vert \mathrm{%
grad\,}\rho \right\vert ^{2}\right) \left\langle X,X\right\rangle
\end{equation*}%
is negative and uniformly bounded away from $0$ along the given horizontal
geodesic.
\end{proposition}

To control the final vertizontal curvature term $\left\langle \nabla
_{X}\left( T_{\frac{\mathrm{grad\,}\rho }{\left\vert \mathrm{grad\,}\rho
\right\vert }}\frac{\mathrm{grad\,}\rho }{\left\vert \mathrm{grad\,}\rho
\right\vert }\right) ,X\right\rangle $ we use compactness to get

\begin{lemma}
\label{L22} For any horizontal lift $\gamma :\left( -\infty ,\infty \right)
\longrightarrow M$ of a geodesic in $B$ that passes through regular points
of $\rho $ and for any $\varepsilon >0,$ there is an interval $\left[ a,b%
\right] $ so that 
\begin{equation*}
\left\vert \frac{1}{b-a}\int_{a}^{b}\left\langle \nabla _{X}\left( T_{\frac{%
\mathrm{grad\,}\rho }{\left\vert \mathrm{grad\,}\rho \right\vert }}\frac{%
\mathrm{grad\,}\rho }{\left\vert \mathrm{grad\,}\rho \right\vert }\right)
,X\right\rangle \right\vert <\varepsilon .
\end{equation*}
\end{lemma}

\begin{proof}
We write 
\begin{eqnarray*}
\left\langle \nabla _{X}\left( T_{\frac{\mathrm{grad\,}\rho }{\left\vert 
\mathrm{grad\,}\rho \right\vert }}\frac{\mathrm{grad\,}\rho }{\left\vert 
\mathrm{grad\,}\rho \right\vert }\right) ,X\right\rangle
&=&D_{X}\left\langle T_{\frac{\mathrm{grad\,}\rho }{\left\vert \mathrm{grad\,%
}\rho \right\vert }}\frac{\mathrm{grad\,}\rho }{\left\vert \mathrm{grad\,}%
\rho \right\vert },X\right\rangle -\left\langle T_{\frac{\mathrm{grad\,}\rho 
}{\left\vert \mathrm{grad\,}\rho \right\vert }}\frac{\mathrm{grad\,}\rho }{%
\left\vert \mathrm{grad\,}\rho \right\vert },\nabla _{X}X\right\rangle \\
&=&D_{X}\left\langle T_{\frac{\mathrm{grad\,}\rho }{\left\vert \mathrm{grad\,%
}\rho \right\vert }}\frac{\mathrm{grad\,}\rho }{\left\vert \mathrm{grad\,}%
\rho \right\vert },X\right\rangle .
\end{eqnarray*}%
Here 
\begin{equation*}
\left\langle T_{\frac{\mathrm{grad\,}\rho }{\left\vert \mathrm{grad\,}\rho
\right\vert }}\frac{\mathrm{grad\,}\rho }{\left\vert \mathrm{grad\,}\rho
\right\vert },\nabla _{X}X\right\rangle =0
\end{equation*}%
since $T_{\frac{\mathrm{grad\,}\rho }{\left\vert \mathrm{grad\,}\rho
\right\vert }}\frac{\mathrm{grad\,}\rho }{\left\vert \mathrm{grad\,}\rho
\right\vert }$ is horizontal and $\nabla _{X}X$ is proportional to $\mathrm{%
grad\,}\rho .$

If $\gamma $ were periodic with period $p,$ then 
\begin{eqnarray*}
\frac{1}{p}\int_{0}^{p}D_{X}\left\langle T_{\frac{\mathrm{grad\,}\rho }{%
\left\vert \mathrm{grad\,}\rho \right\vert }}\frac{\mathrm{grad\,}\rho }{%
\left\vert \mathrm{grad\,}\rho \right\vert },X\right\rangle &=&\frac{1}{p}%
\left. \left\langle T_{\frac{\mathrm{grad\,}\rho }{\left\vert \mathrm{grad\,}%
\rho \right\vert }}\frac{\mathrm{grad\,}\rho }{\left\vert \mathrm{grad\,}%
\rho \right\vert },X\right\rangle \right\vert _{0}^{p} \\
&=&0.
\end{eqnarray*}%
By compactness of the unit tangent bundle of $M,$ we can find an interval on
which $\gamma $ is almost periodic with the error being any predetermined
quantity, yielding the result.
\end{proof}

Combining the previous three results with our vertizontal curvature formula (%
\ref{AF}), we get\newline

\begin{proposition}
Let $\pi :M\longrightarrow B$ be a horizontally homothetic submersion of a
compact Riemannian manifold with nonconstant dilation function. Any
horizontal lift of a geodesic in $B$ that passes through regular points and
is sufficiently close to the minimum of $\rho $ passes through some points
with some negative vertizontal sectional curvatures.
\end{proposition}

Our main theorem is a corollary of the previous proposition.\medskip

\section*{Nonpositive Curvature}

Much of our argument carries through to compact manifolds with nonpositive
sectional curvature; however, the two vertizontal terms 
\begin{equation}
-2\left\langle T_{\frac{\mathrm{grad\,}\rho }{\left\vert \mathrm{grad\,}\rho
\right\vert }}\left( \nabla _{X}\frac{\mathrm{grad\,}\rho }{\left\vert 
\mathrm{grad\,}\rho \right\vert }\right) ^{vert},X\right\rangle
-\left\langle T_{\frac{\mathrm{grad\,}\rho }{\left\vert \mathrm{grad\,}\rho
\right\vert }}X,T_{\frac{\mathrm{grad\,}\rho }{\left\vert \mathrm{grad\,}%
\rho \right\vert }}X\right\rangle  \notag
\end{equation}%
are nonpositive and one expects them to usually be negative. By assuming
them away we get

\begin{theorem}
Every horizontally homothetic submersion of a compact manifold with
nonpositive sectional curvature is a Riemannian submersionn (up to a change
of scale on the base), provided the fibers are totally geodesic.
\end{theorem}

\begin{proof}
Since the fibers are totally geodesic and $\frac{\mathrm{grad\,}\rho }{%
\left\vert \mathrm{grad\,}\rho \right\vert }$ is in the kernel of the $A$%
--tensor%
\begin{equation}
\left\langle R\left( X,\frac{\mathrm{grad\,}\rho }{\left\vert \mathrm{grad\,}%
\rho \right\vert }\right) \frac{\mathrm{grad\,}\rho }{\left\vert \mathrm{%
grad\,}\rho \right\vert },X\right\rangle =-\left( D_{\frac{\mathrm{grad\,}%
\rho }{\left\vert \mathrm{grad\,}\rho \right\vert }}\left\vert \mathrm{grad\,%
}\rho \right\vert +\left\vert \mathrm{grad\,}\rho \right\vert ^{2}\right)
\left\langle X,X\right\rangle .  \notag
\end{equation}%
If $\rho $ is nonconstant and we are close enough to a point where it
achieves its maximum, then 
\begin{eqnarray*}
\left\langle R\left( X,\frac{\mathrm{grad\,}\rho }{\left\vert \mathrm{grad\,}%
\rho \right\vert }\right) \frac{\mathrm{grad\,}\rho }{\left\vert \mathrm{%
grad\,}\rho \right\vert },X\right\rangle &=&-\left( D_{\frac{\mathrm{grad\,}%
\rho }{\left\vert \mathrm{grad\,}\rho \right\vert }}\left\vert \mathrm{grad\,%
}\rho \right\vert +\left\vert \mathrm{grad\,}\rho \right\vert ^{2}\right)
\left\langle X,X\right\rangle \\
&>&0.
\end{eqnarray*}%
So $\rho $ must be constant if the domain is compact, the sectional
curvature is nonpositive, and the fibers are totally geodesic.
\end{proof}

\textbf{Acknowledgement: }We are grateful to Paula Bergen for copyediting a
draft of this \medskip paper.

\end{document}